\documentclass{article}

\usepackage{amsfonts}
\usepackage{amssymb}
\usepackage{textcomp}
\usepackage{latexsym,amsmath}
\usepackage{graphics}

\usepackage{amsfonts}
\usepackage{amssymb}
\usepackage{textcomp}
\usepackage{latexsym,amsmath}

\renewcommand{\phi}{\varphi}
\newcommand{\be}{\begin{equation}}
\newcommand{\ee}{\end{equation}}
\newcommand{\ba}{\begin{eqnarray}}
\newcommand{\ea}{\end{eqnarray}}
\newcommand{\ban}{\begin{eqnarray*}}
\newcommand{\ean}{\end{eqnarray*}}

\newcommand{\nul}{{\bf0}}
\newcommand{\rd}{{\mathbb R}^d}
\newcommand{\zd}{{\mathbb Z}^d}
\newcommand{\td}{{\mathbb T}^d}

\newcommand{\z} {{\mathbb Z}}

\newcommand{\n} {{\mathbb N}}

\newcommand{\h}{\widehat}
\newcommand{\w}{\widetilde}

\def\supp{\operatorname{supp}}
\def\sinc{\operatorname{sinc}}
\def\const{\operatorname{const}}

\def\N{{{\Bbb N}}}
\def\Z{{{\Bbb Z}}}
\def\T{{{\Bbb T}}}
\def\R{{\Bbb R}}

\def\vp{{\varphi}}

\def\({\left(}
\def\){\right)}

\def\l{{\lambda }}
\def\a{{\alpha }}

\def\a{{\alpha}}
\def\b{{\beta}}
\def\d{{\delta}}

\def\vp{{\varphi}}

\def\g{{\gamma }}

\def\F{{\mathcal{F}}}
\def\A{{\mathcal{A}}}
\def\S{{\mathcal{S}}}

\def\B{{\mathcal{B}}}

\newtheorem{theo}{Theorem}
\newtheorem{lem}[theo]{Lemma}

\newtheorem {coro} [theo] {Corollary}
\newtheorem {defi} [theo] {Definition}
\newtheorem {rem} [theo] {Remark}

\title{Uniform approximation by  multivariate \\ quasi-projection operators
}
\author{
Yu. Kolomoitsev$^{1, 2}$ and M. Skopina$^{3, 4}$
}
\date{\small $^{1}$Universit\"at zu L\"ubeck,
Institut f\"ur Mathematik, L\"ubeck, Germany; kolomoitsev@math.uni-luebeck.de  \\
\small $^{2}$Institute of Applied Mathematics and Mechanics of NAS of Ukraine, Slov'yans'k, Ukraine
\\
$^{3}$St. Petersburg State University,  Saint Petersburg,  Russia;
skopina@ms1167.spb.edu
\\
$^{4}$Regional Mathematical Center of Southern Federal University
}

\begin{document}

\maketitle

\begin{abstract}
Approximation properties   of quasi-projection operators  $Q_j(f,\varphi, \widetilde{\varphi})$ are studied. Such an operator
is associated with a function $\varphi$ satisfying the Strang-Fix conditions and a  tempered distribution $\widetilde{\varphi}$
such that compatibility conditions with $\varphi$ hold.  Error estimates in the uniform norm  are obtained for a wide class of quasi-projection operators defined on  the space of uniformly continuous functions and on the anisotropic Besov spaces.
Under additional assumptions on $\varphi$ and $\widetilde{\varphi}$, two-sided estimates in terms of realizations of the $K$-functional are  also obtained.
\end{abstract}

\bigskip


\medskip

	
\section{Introduction}

Quasi-projection operators are a generalisation of so-called scaling expansions
$$
(Q_jf)(x)=2^j\sum_{k\in\z} \langle f, {\w\phi}(2^j\cdot-k)\rangle \phi(2^jx-k) , \quad f, \phi, \w\phi\in L_2({\mathbb R}),
$$
playing an important role in the wavelet theory (see, e.g.,~\cite{BDR, Jia2, KPS, KS, Sk1}).  Such expansions are
also well defined whenever the inner product $\langle f, {\w\phi}(2^j\cdot-k)\rangle$   has meaning and the series
converges in some sense. In~\cite{Jia2}, Jia considered a larger class of quasi-projection operators $Q_j$ and obtained error estimates in $L_p$ and other function spaces for $Q_j$.
The classical Kantorovich-Kotelnikov operators (see, e.g.,~\cite{ CV0, CV2, KS3, KS20, OT15, VZ2})  are operators of the same form  $Q_j$  with  the characteristic function of $[0,1]$ as $\w\phi$. Another classical special case of quasi-projection operators
is the sampling expansion
$$
(Q_wf)(x)=\sum_{k\in\z} f(w^{-1}k)\phi(w x-k)=\sum_{k\in\z}  \langle f, {\delta}(w\cdot-k)\rangle \phi(w x-k),
$$
where  $\delta$ is the Dirac delta-function, $\phi(x)=\sinc x:=\sin \pi x/\pi x$, and $w>0$. Since $\delta$ is a tempered distribution, under the usual notation
$\langle f, {\delta}\rangle:=\delta(f)$, the operator $Q_w$ is defined only for functions $f$ from the Schwartz class,
but to extend this class, one can set $\langle f, {\delta}\rangle:=\langle \h f, \h\delta\rangle$.
The sampling expansion  is of great applied importance, especially actively it is  used by
engineers working in  signal processing.  The sampling-type operators $Q_w$
associated with different functions $\phi$ and their approximation properties as $w\to\infty$  were studied by a
lot of authors (see, e.g.,~\cite{Butz4, Butz6, Brown,  BD,   Butz5,   Butz7, KM, {NU}, SS, Unser}).

Given a matrix $M$, we define the multivariate quasi-projection operator $Q_j(f, \phi,\w\phi)$ associated with
a function $\phi$ and a distribution/function $\w\phi$ as follows
$$
 Q_j(f, \phi,\w\phi)(x)=|\det M|^{j}  \sum_{n\in\zd} \langle f,\w\vp(M^j\cdot-n)\rangle \vp(M^j x-n),
$$
where the "inner product" $\langle f,\w\vp(M^j\cdot-n) \rangle$ has meaning in some sense.
If the Fourier transform of $f$ has enough decay, then  the  operators $Q_j(f, \phi,\w\phi)$
with $\langle f,\w\vp\rangle:=\langle \h f,\h{\w\vp\rangle}$ are well defined
for a wide class of tempered distributions $\w\phi$ and appropriate functions $\phi$.
Approximation properties of such quasi-projection operators were studied for fast decaying  functions $\phi$ in~\cite{KS}.
The error estimates in the $L_p$-norm, $2\le p\le\infty$, were given in terms of the Fourier transform of $f$,
and the approximation order of $Q_j(f, \phi,\w\phi)$  was found for the isotropic matrices $M$.
  Similar results for a class of bandlimited $\phi$ and $p<\infty$ were obtained in~\cite{KKS}, and then the
  class of functions $\phi$ was  essentially extended  in~\cite{cksv}. These results were improved  in several directions
in~\cite{KS19}  and~\cite{KS20+} respectively for  fast decaying and bandlimited $\phi$.
The  error estimates  in $L_p$-norm,   $1\le p<\infty$,  were given  in  terms of  moduli of smoothness and best approximations,
and the class of  approximated function $f$ was extended.
Also the requirement on smoothness of $\h\phi$ and $\h{\w\phi}$ was  weakened in~\cite{KS20+} due to
using the  Fourier multipliers method.

The technique developed in~\cite{KS19}, \cite{KS20+}  does not work for $p=\infty$ in full because a function from $L_\infty$
cannot be uniformly approximated by functions from $L_2$.
 In~\cite{KS20+}, the case $p=\infty$ is  considered,
but the estimates are obtained only for $f$ satisfying the additional assumption $f(x)\to0$ as $|x|\to \infty$.
The goal of the present paper is to fix this drawback.  It succeeded due to the using a new technique
based on convolution representations. In particular,  the "inner product" $\langle f,\w\vp(M^j\cdot+n) \rangle$
is defined as a limit of some convolutions, so that the operator  $Q_j(f, \phi,\w\phi)$  is well defined on the space of uniformly
continuous functions and on the anisotropic Besov spaces whenever the series $\sum_{k\in\zd}|\phi(x-k)|$ converges uniformly on any compact.
Error estimates in the uniform norm for  $Q_j(f, \phi,\w\phi)$ are obtained
under the assumptions  of the Strang-Fix conditions for $\phi$, compatibility conditions of $\phi$ with $\w\phi$,
and the belonging of particular functions related to $\phi$ and $\w\phi$ to Wiener's algebra. The results are very similar
to those for $p<\infty$ in~\cite{KS20+}, where some Fourier multiplier conditions
are assumed instead of belonging to Wiener's algebra. Under additional assumptions on $\phi$ and $\w\phi$, two-sided estimates in terms of
realizations of the $K$-functional are presented. A new Whittaker--Nyquist--Kotelnikov--Shannon-type theorem is also proved.

The paper is organized as follows. Notation  and preliminary information are given in Sections~2 and 3, respectively.
Section~4 contains auxiliary results.  The main results are presented in Section~5. The case of weak compatibility of $\phi$ and $\w\phi$ is discussed in Subsection~5.2. Subsection~5.3 is devoted to approximation by operators $Q_j(f,\phi,\w\phi)$ in the case of strict compatibility of $\phi$ and $\w\phi$.

\section{Notation}

We use the standard multi-index notations.
    Let $\n$ be the set of positive integers, $\rd$ be the $d$-dimensional Euclidean space,
    $\zd$ is the integer lattice  in $\rd$,
    $\td=\rd\slash\zd$ be the $d$-dimensional torus.
    Let  $x = (x_1,\dots, x_d)^{T}$ and
    $y =(y_1,\dots, y_d)^{T}$ be column vectors in $\rd$,
    then $(x, y):=x_1y_1+\dots+x_dy_d$,
    $|x| := \sqrt {(x, x)}$; $\nul=(0,\dots, 0)^T\in \rd$;  		
		$\z_+^d:=\{x\in\zd:~x_k\geq~{0}, k=1,\dots, d\}.$
		If  $r>0$, then $B_r$
		denotes the ball of radius $r$ with the center in $\nul$.

If $\alpha\in\zd_+$, $a,b\in\rd$, we set
    $$
    [\alpha]=\sum\limits_{j=1}^d \alpha_j, \quad
    D^{\alpha}f=\frac{\partial^{[\alpha]} f}{\partial x^{\alpha}}=\frac{\partial^{[\alpha]} f}{\partial^{\alpha_1}x_1\dots
    \partial^{\alpha_d}x_d}, \,\,\quad
 a^b=\prod\limits_{j=1}^d a_j^{b_j},\,\,\quad
  \alpha!=\prod\limits_{j=1}^d\alpha_j!.
    $$

If $M$ is a $d\times d$ matrix,
then $\|M\|$ denotes its operator norm in $\rd$; $M^*$ denotes the conjugate matrix to $M$, $m=|\det M|$;
the identity matrix is denoted by $I$.

A $d\times d$ matrix $M$ whose
eigenvalues are bigger than $1$ in modulus is called a  dilation matrix.
 We denote the set of all dilation matrices by $\mathfrak{M}$.
It is well known that $\lim_{j\to\infty}\|M^{-j}\|=0$ for any dilation matrix $M$.

A  matrix $M$ is called isotropic if it is similar to a diagonal matrix
such that numbers $\lambda_1,\dots,\lambda_d$ are placed on the main diagonal
and $|\lambda_1|=\cdots=|\lambda_d|$.

By $L_p$ we  denote the space $L_p(\rd)$, $1\le p\le \infty$, with the norm $\Vert \cdot\Vert_p=\Vert \cdot\Vert_{L_p(\R^d)}$.
As usual, ${C}$ denotes the space of all uniformly continuous bounded functions on $\R^d$  equipped with the norm $\Vert f\Vert=\max_{x\in \R^d}|f(x)|$.

We use $W_p^n$, $1\le p\le\infty$, $n\in\n$, to denote
 the Sobolev space on~$\rd$, i.e. the set of
functions whose derivatives up to order $n$ are in $L_p$, with usual Sobolev norm.

If $f, g$ are functions defined on $\rd$ and $f\overline g\in L_1$,
then
$$
\langle  f, g\rangle:=\int_{\rd}f(x)\overline{g(x)}dx.
$$
As usual, the convolution for appropriate functions $f$ and $g$ is defined by
$$
f*g(x)=\int_{\R^d} f(t)g(x-t)dt.
$$

If $f\in L_1$,  then its Fourier transform is
$$
\mathcal{F}f(\xi)=\widehat
f(\xi)=\int_{\rd} f(x)e^{-2\pi i
(x,\xi)}\,dx.
$$

For any function $f$, we denote $f^-(x)=\overline{f(-x)}$ .

Denote by $\mathcal{S}$ the Schwartz class of functions defined on $\rd$.
    The dual space of $\mathcal{S}$ is $\mathcal{S}'$, i.e. $\mathcal{S}'$ is
    the space of tempered distributions.
       Suppose $f\in \mathcal{S}$ and $\phi \in \mathcal{S}'$, then
  $\langle f, \phi \rangle:=\phi(f)$.
    For any  $\phi\in \mathcal{S}'$, we define $\overline{\vp}$ and $\vp^-$  by $\langle f, \overline{\phi}\rangle:=\overline{\langle \overline{f},  \phi\rangle}$,
    $f\in \mathcal{S}$,
    and $\langle f, {\phi}^-\rangle:=\overline{\langle {f}^-,  \phi\rangle}$,
    $f\in \mathcal{S}$, respectively. The Fourier transform of $\phi$ is
    defined by $\langle \h f, \h \phi\rangle=\langle f, \phi\rangle$,
    $f\in \mathcal{S}$. The convolution of $\phi \in \mathcal{S}'$ and $f\in \mathcal{S}$ is given by
$f*\vp:=\mathcal{F}^{-1}(\h\phi\h f)$. 

We say that  $f\in \mathcal{S}'$ belongs to Wiener's algebra $W_0$, if there exists a function $g\in L_1$ such that
\begin{eqnarray}\label{belw}
f(x)=\int_{\R^d}g(\xi)e^{i(x,\xi)}d\xi.
\end{eqnarray}
The corresponding norm is given by  $\|f\|_{W_0}=\|g\|_{1}$.
		
 For a fixed matrix $M\in \mathfrak{M}$ and  a function $\phi$  defined on $\rd$, we set
$$
\phi_{jk}(x):=m^{j/2}\phi(M^jx+k),\quad j\in\z,\quad k\in\rd.
$$
If  $\w\phi\in \mathcal{S}'$, $j\in\z, k\in\zd$, then we define $\w\phi_{jk}$ by
        $$
        \langle f, \w\phi_{jk}\rangle:=
        \langle f_{-j,-M^{-j}k},\w\phi\rangle,\quad  f\in \mathcal{S}.
        $$

Denote by ${\cal L}_\infty$ the set of functions $\phi\in L_\infty$
such that $\sum_{k\in\zd} |\phi(\cdot+k)|\in L_\infty(\td)$, and by ${\cal L}C$ the set of functions $\vp\in C(\R^d)$ such that
the series $\sum_{k\in\zd} \left|\phi(\cdot+k)\right|$ converges uniformly on any compact set. Both ${\cal L}_\infty$
and ${\cal L}C$ are Banach spaces with the norm
$$
\|\phi\|_{{\cal L}C}=\|\phi\|_{{\cal L}_\infty}:=	\bigg\|\sum_{k\in\zd} \left|\phi(\cdot+k)\right|\bigg\|_{L_\infty(\td)}.
$$

For any $d\times d$ matrix $A$, we introduce the space
$$
\mathcal{B}_{A}:=\{g\in L_\infty\,:\, \supp \h g\subset A^*\T^d \}
$$
and the corresponding anisotropic best approximations
$$
E_{A} (f):=\inf\{\Vert f-g\Vert\,:\, g\in \mathcal{B}_{A}\}.
$$

Let $\alpha$ be a positive function defined on the set of all $d\times d$ matrices $A$.
We consider the following anisotropic Besov-type space associated with a matrix~$A$.
  We  say that
$f\in \mathbb{B}_{A}^{\a(\cdot)}$ if
$$
\Vert f\Vert_{\mathbb{B}_{A}^{\a(\cdot)}}:=\Vert f\Vert+\sum_{\nu=1}^\infty\a(A^\nu) E_{A^\nu} (f)<\infty.
$$
Note that in the case $A=2 I_d$ and $\a(\cdot)\equiv \a_0\in \R$, the space $\mathbb{B}_{A}^{\a(\cdot)}$ coincides with the classical Besov space $B_{\infty,1}^{\a_0}(\R^d)$.

For any matrix $M\in \mathfrak{M}$, we denote by  $\mathcal{A}_M$ the set of all positive functions $\a\,:\,\R^{d\times d}\to \R_+$ that satisfy the condition $\a(M^{\mu+1})\le c(M)\a(M^\mu)$ for all $\mu\in \Z_+$.

For any $d\times d$ matrix $A$, we introduce also the anisotropic fractional modulus of smoothness of order $s$, $s>0$,
$$
\Omega_s(f,A):=\sup_{|A^{-1}t|<1, t\in \R^d} \Vert \Delta_t^s f\Vert,
$$
where
$$
\Delta_t^s f(x):=\sum_{\nu=0}^\infty (-1)^\nu \binom{s}{\nu} f(x+t\nu).
$$
Recall that the standard fractional modulus of smoothness of order $s$, $s>0$, is defined
by
\begin{equation}\label{fm}
 \omega_s(f,h):=\sup_{|t|<h} \Vert \Delta_t^s f\Vert,\quad h>0.
\end{equation}
We refer to~\cite{KT20} for the collection of basic properties of moduli of smoothness in $L_p(\R^d)$.

For an appropriate function $f$ and  $s>0$, the fractional power of Laplacian is given by
$$
(-\Delta)^{s/2} f(x):= \mathcal{F}^{-1}\(|\xi|^s\h f(\xi)\)(x).
$$

As usual, if $\{a_k\}_{k\in \Z^d}$ is a sequence, then $\left\Vert \{a_k\}_{k}\right\Vert_{\ell_{\infty}}:=\sup\limits_{k\in \Z^d}|a_k|$.

By $\eta$ we denote a real-valued function in $C^\infty(\R^d)$ such that $\eta(\xi)=1$ for $\xi \in \T^d$ and $\eta(\xi)=0$ for $\xi\not\in 2\T^d$. Next, for $\d>0$  and a $d\times d$ matrix $A$, we set
$$
\eta_\d=\eta(\d^{-1}\cdot)\quad\text{and}\quad \mathcal{N}_\d=\mathcal{F}^{-1}\eta_\d,
$$
$$
\eta_A=\eta(A^{*-1}\cdot) \quad\text{and}\quad \mathcal{N}_A=\mathcal{F}^{-1}\eta_A.
$$

\section{ Preliminary information and main definitions.}
\label{sa}

In what follows, we discuss the quasi-projection operators
$$
Q_j(f,\phi,\w\phi):=\sum_{k\in\zd} \langle f, {\w\phi}_{jk}\rangle \phi_{jk},
$$
where the "inner product" $\langle f, {\w\phi}_{jk}\rangle$ is defined in a special way and the series converges in some sense.

The expansions $\sum_{k\in\zd} \langle f, {\w\phi}_{jk}\rangle \phi_{jk}$ are elements of the shift-invariant
spaces generated by~$\phi$. It is well known that a function $f$ can be approximated by elements of such  spaces only if $\phi$
satisfies a special property, the so-called Strang-Fix conditions.

\begin{defi}\label{d2}
A  function $\phi$ is said to satisfy {\em the Strang-Fix conditions} of
order $s$  if $D^{\beta}\h{\phi}(k) = 0$ for every $\beta\in\z_+^d$, $[\beta]<s$,
and for all $k\in\zd\setminus \{\nul\}$.
\end{defi}

Certain compatibility conditions for a distribution $\w\phi$ and a function $\phi$
is also required to provide good approximation properties of the operator $Q_j(f,\phi,\w\phi)$. For our purposes, we will use the following conditions.

\begin{defi}
\label{d3}
 A tempered distribution  $\w\phi$ and a function $\phi$ is said to be  {\em weakly compatible of order~$s$} if $D^{\beta}(1-\h\phi\h{\w\phi})({\bf 0}) = 0$ 	for every $\beta\in\z_+^d$, $[\beta]<s$.
\end{defi}

\begin{defi}
\label{d1}
 A tempered distribution  $\w\phi$ and a function $\phi$ is said to be  {\em strictly compatible} if there exists $\delta>0$ such that  $\overline{\h\phi}(\xi)\h{\w\phi}(\xi)=1$  a.e. on $\delta\td$.
\end{defi}

Denote by $\mathcal{S}_N'$ the set of  tempered distribution $\w\phi$
	whose Fourier transform $\h{\w\phi}$ is a measurable function on $\rd$
	such that $|\h{\w\phi}(\xi)|\le C_{\w\phi}(1+ |\xi|)^{N}$
	 for almost all $\xi\in\rd$, $N\ge 0$.
For $\w\phi\in \mathcal{S}'_N$ and appropriate classes
of  functions $\phi$,  the quasi-projection operators $Q_j(f,\phi,\w\phi)$ with
$\langle f, {\w\phi}_{jk}\rangle:=\langle \h f, \h{\w\phi_{jk}}\rangle$
were studied  in~\cite{KKS},  \cite{KS20}, \cite{KS}, and~\cite{Sk1}.
In particular,  approximation by such operators in the uniform norm was considered in~\cite{KS}. The following statement can be derived from Theorems~4 and~5 in~\cite{KS}.

\medskip

\noindent{\bf Theorem A.}	
	{\it Let $s\in\n$,  $N\ge0$,  $\delta\in(0, 1/2)$ and  $M\in  \mathfrak{M}$.
	Suppose
\begin{itemize}
 	\item[$1)$]  $\phi, \h\phi \in {\cal L}_\infty$;

      \item[$2)$]  $\h\phi(\cdot+l)\in C^{s}(B_\delta)$  for all $l\in\zd\setminus\{\nul\}$
and $\sum\limits_{l\ne\nul}\sum\limits_{\|\beta\|_1=s}
		\sup\limits_{|\xi|<\delta}|D^\beta\h\phi(\xi+l)|< \infty$;

	  \item[$3)$]  the Strang-Fix conditions of order $s$ are satisfied	for $\phi$;

     \item[$4)$] $\w\phi \in \mathcal{S}_N'$ and $\overline{\h\phi}\h{\w\phi}\in C^s(B_\delta)$;
 	
   \item[$5)$]  $\phi$ and ${\w\phi}$ are weakly  compatible of order~$s$.
\end{itemize}

\noindent
 		If $ f\in L_{\infty}$ is such that $\h f\in L_1$ and
$\h f(\xi)=\mathcal{O}(|\xi|^{-N-d-\varepsilon})$, $\varepsilon>0$,
as $|\xi|\to\infty$,  then
$$
\bigg\| f- \sum_{k\in\zd}\langle \h f, \h{\w\phi_{jk}}\rangle\phi_{jk}\bigg\|_{\infty} \le C_1\|M^{*-j}\|^N  \int\limits_{|M^{*-j}\xi|\ge\delta}|\xi|^N |\h f(\xi)|\,d\xi
+ C_2 \|M^{*-j}\|^s  \int\limits_{|M^{*-j}\xi|\le\delta}|\xi|^s |\h f(\xi)|\,d\xi,
$$
where the constants $C_1$ and $C_2$ do not depend on $f$ and $j$.
If, moreover,  $M$ is an isotropic matrix, then
	 $$
\bigg\| f- \sum_{k\in\zd}\langle \h f, \h{\w\phi_{jk}}\rangle\phi_{jk}\bigg\|_{\infty} \le C_3
	 \begin{cases}
	 |\lambda|^{-j(N + \varepsilon)}  &\mbox{if }
	s> N + \varepsilon\\
	  (j+1) |\lambda|^{-js} &\mbox{if }
	 s= N+ \varepsilon \\
	|\lambda|^{-js}
	 &\mbox{if }
	 s< N+ \varepsilon
	\end{cases},
$$
	where $\lambda$ is an eigenvalue of $M$ and  $C_3$ does not depend on $j$.
}

\medskip

Unfortunately, there is a restriction on  the decay of  $\h f$ in this theorem.
 Obviously, such a restriction is redundant for some special cases, in particular,
 the inner product $\langle f, {\w\phi}_{jk}\rangle$ has meaning for any $f\in L_\infty$ whenever  $\w\phi $ is an integrable function. Moreover,  Theorem~A provides approximation order for $Q_j(f,\phi,\w\phi)$ only for isotropic matrices $M$, and even for this case,  more accurate  error estimates in terms of smoothness of $f$ were not obtained. The mentioned drawbacks  were avoided
in~\cite[Theorem~17']{KS3}, where  the uniform approximation by quasi-projection  operators $Q_j(f,\phi,\w\phi)$ associated with a summable function $\w\phi$ and a bandlimited function $\phi$ was investigated.
 To formulate this result, we need to introduce the space  ${\cal B}$ consisting
of  functions $\phi$ such that
$
\phi=\F^{-1}\theta,
$
where the function   $\theta$ is supported in a rectangle $ R\subset\rd$ and
	$\theta\big|_{R}\in C^d(R)$.

\medskip

\noindent{\bf Theorem B.}
{\it Let  $s\in\n$, $\delta, \delta'>0$,  $M\in  \mathfrak{M}$ and  $f\in L_\infty$.  Suppose
\begin{itemize}
  \item[$1)$] $\phi\in\cal B\cap {\cal L}_\infty$, $\supp\h\phi\subset B_{1-\delta'}$, and $\h\phi\in C^{s+d+1}(B_\delta)$;
 \item[$2)$]  $\w\phi \in L_1$ and  $\h{\w\phi}\in C^{s+d+1}(B_\delta)$;
 \item[$3)$] $\phi$ and ${\w\phi}$ are weakly  compatible of order~$s$.
\end{itemize}
 Then
\begin{equation*}
  \bigg\|f-\sum\limits_{k\in\zd}
\langle f,\widetilde\phi_{jk}\rangle \phi_{jk}\bigg\|_\infty\le C\,\omega_s\(f,\|M^{-j}\|\),
\end{equation*}
where $C$ does not depend on $f$ and $j$.
}

\medskip

In what follows, we will consider  a class of quasi-projection operators $Q_j(f,\vp,\w\vp)$ associated with a tempered distribution $\w\phi$ belonging to the class $\mathcal{S}_{\a;M}'$, where  $M\in \mathfrak{M}$ and $\a\in \mathcal{A}_M$. We say that  $\w\phi\in\mathcal{S}_{\a;M}'$  if $\h{\w\phi}$ is a measurable locally bounded function and, for any  $\nu\in \Z_+$,
the function $\mathcal{N}_{M^{\nu}}*\w\vp^-$ is summable and
\begin{equation}\label{DefS}
  \|\mathcal{N}_{M^{\nu}}*\w\vp^-\|_1\le c\,\a(M^\nu),
\end{equation}
where $c$ is independent of $\nu$.

Let us show that  inequality~\eqref{DefS} is satisfied for the most  important special cases of $\w\vp$.
Since $\|\mathcal{N}_{M^{\nu}}\|_1=\|\mathcal{N}_{1}\|_1
=\|\h\eta\|_1$, one can easily see that~\eqref{DefS} with $\alpha\equiv 1$ holds true if $\w\vp\in L_1$ or $\w\phi$ is the Dirac delta-function $\delta$. Let now $\w\phi$ be a distribution associated with
the differential operator $D^\beta$, $\beta\in \z_+^d$, i.e.,  $\w\vp(x)=(-1)^{[\beta]}D^\b \d(x)$ (see~\cite{KKS}). In this case, we have
$$
\mathcal{N}_{M^{\nu}}*\w\vp^-={\mathcal F}^{-1}(\h{\mathcal{N}_{M^{\nu}}}\h{D^\b \d})=
{\mathcal F}^{-1}(\h{D^\b\mathcal{N}_{M^{\nu}}})=
D^\b\mathcal{N}_{M^{\nu}}.
$$
If  $M={\rm diag}(m_1,\dots,m_d)$ and
$\alpha(M)=m_1^{\beta_1}\dots m_d^{\beta_d}$, then  using
 Bernstein's inequality (see, e.g.,~\cite[p.~252]{Timan}), we have
$$
\|\mathcal{N}_{M^{\nu}}*\w\vp^- \|_1=\| D^\b\mathcal{N}_{M^{\nu}}\|_1\le m_1^{\beta_1}\dots m_d^{\beta_d}\|\mathcal{N}_{M^{\nu}}\|_1=\alpha(M)\|\h\eta\|_1.
$$
Similarly,  if $M$ is an isotropic matrix, then $\w\phi$ belongs to the class $\mathcal{S}_{\a;M}'$ with  $\alpha(M)=m^{{[\beta]}/d}$.

To extend the operator $Q_j(f,\phi,\w\phi)$ associated  with $\w\phi \in \mathcal{S}_{\a;M}'$
onto the Besov spaces $\mathbb{B}_{M}^{\a(\cdot)}$ and onto~${C}$, we need to
define the "inner product" $\langle f,\w\phi_{jk}\rangle$ properly.
A similar extension for the case $p<\infty$ was realized in~\cite{KS19}  and~\cite{KS20+},
but the definition of $\langle f,\w\phi_{jk}\rangle$ given there is not appropriate for us now.
	
\begin{defi}
\label{def0}
	Let $M\in \mathfrak{M}$, $\d\in(0,1]$,
	$\w\vp\in \mathcal{S}_{\a;M}'$ with $\a\in \mathcal{A}_M$, and the functions $\{T_\mu\}_{\mu\in\Z_+}$ be such that $T_\mu\in \mathcal{B}_{\d M^\mu}$ and
\begin{equation}\label{eT}
  \|f-T_\mu\|_{}\le c(d){E}_{\d M^\mu}(f)_{}.
\end{equation}	
For every $f\in \mathbb{B}_{M}^{\a(\cdot)}$ and for every $f\in C$ in the case $\a\equiv {\const}$,  we set
\begin{equation}\label{con}
  \langle f,\w\vp_{0k}\rangle:
  =\lim_{\mu\to\infty} T_\mu*(\mathcal{N}_{M^{\mu}}*\w\vp^-)(k), \quad k\in\zd,
\end{equation}
and
	$$
	\langle f,\w\vp_{jk}\rangle:= m^{-j/2}\langle f(M^{-j}\cdot),\w\vp_{0k}\rangle, \quad j\in\z_+.
	$$
	\end{defi}

To approve this definition, we note that, according to Lemma~\ref{lem1} (which will be presented latter), the limit in~\eqref{con} exists and does not depend on a choice of $\d$ and functions $T_\mu$. 

\begin{rem}
\label{prop003}
 Note that since $\mathcal{N}_{M^{\mu}}*\w\vp^-\in L_2$, the convolution in~\eqref{con} is also well-defined for all $T_\mu\in L_2$ and
\begin{equation}\label{conv}
  T_\mu*(\mathcal{N}_{M^{\mu}}*\w\vp^-)(k)=\langle \h{T_\mu},\h{\w\vp_{0,-k}}\rangle.
\end{equation}
Further, if a function $f$ satisfies the conditions of Theorem~A,  then the functional $\langle f,\w\vp_{0k}\rangle$  defined in the sense of
Definition~\ref{def0} is the same as in that theorem. Indeed, since $\h f\in L_1$,
we have $\lim_{|x|\to\infty} f(x)=0$, and, due to Lemma~15 in~\cite{KS20+}, one can choose functions
$T_\mu\in \mathcal{B}_{\d M^\mu}\cap L_2$ satisfying~\eqref{eT}. Then repeating the arguments of Remark~12 in~\cite{KS19},
 we obtain
$$
  \langle \h f, \h{\w\vp_{0k}}\rangle =\lim_{\mu\to\infty}\langle \h{T_\mu},\h{\w\vp_{0,-k}}\rangle,
$$
which, together with~\eqref{conv}, yields the equality $\langle f,\w\vp_{0k}\rangle=\langle \h f,\h{\w\vp_{0,-k}}\rangle$.
\end{rem}

The main results of this paper will be given in terms of Wiener's algebra $W_0$.
Various conditions of belonging to Wiener's algebra are overviewed in detail in the survey \cite{LST} (see also~\cite{Kol}, \cite{KL13}, \cite{KL17}, for some new efficient sufficient conditions).
Here, we mention only the Beurling-type theorem, which states that if $h\in W_2^{k}$ with $k>d/2$, then $h\in W_0$ (see, e.g.,~\cite[Theorem~6.1]{LST}).

\section{Auxiliary results}
	
\begin{lem}
\label{lem1}
Let  $M\in \mathfrak{M}$, $n\in\N$, $\delta\in (0,1]$, and $\a\in \A_M$.
Suppose that $\w\vp$, $f$, and the functions $T_\mu$, $\mu\in\Z_+$, are as in Definition~\ref{def0}
and
$$
q_\mu(k):=T_\mu*(\mathcal{N}_{M^\mu}*\w\vp^-)(k).
$$
Then the sequence $\{\{q_\mu(k)\}_k\}_{\mu=1}^\infty$ converges in $\ell_\infty$ as $\mu\to\infty$ and its limit does not depend on the choice of $T_\mu $ and $\d$; a fortiori for every $k\in\zd$ there exists a limit
$\lim_{\mu\to\infty}q_\mu(k)$  independent on the choice of $T_\mu $ and $\d$.  Moreover, for all $f\in \mathbb{B}_{M}^{\a(\cdot)}$, we have
\begin{equation}\label{elem1}
	\sum_{\mu=n}^\infty
	\|\{q_{\mu+1}(k)-q_\mu(k)\}_{k}\|_{\ell_\infty}\le
	 c\sum_{\mu=n}^\infty \alpha (M^\mu)
E_{\delta M^\mu}(f)_{},
\end{equation}
where $c$ depends only on $d$ and $M$.
\end{lem}
{\bf Proof.} Denote $\mu_0=\min\{\mu\in \N\,:\,\T^d\subset \frac12 M^{*\nu}\T^d\quad\text{for all}\quad\nu\ge \mu-1\}$.
It is easy to see that $\eta({M^{*-\mu}\cdot})\eta(M^{*-\mu-\mu_0}\cdot)=
\eta(M^{*-\mu}\cdot)$, and hence
$\mathcal{N}_{M^\mu}=\mathcal{N}_{M^\mu} *\mathcal{N}_{M^{\mu+\mu_0}}$. Similarly,
$\mathcal{N}_{M^{\mu+1}}=\mathcal{N}_{M^{\mu+1}} *\mathcal{N}_{M^{\mu+\mu_0}}$,
and taking into account that $T_\mu=T_\mu*N_{M^{\mu}}$ and
$T_{\mu+1}=T_{\mu+1}*\mathcal{N}_{M^{\mu+1}}$, we have
$$
T_\nu*\mathcal{N}_{M^{\nu}} *\w\phi^-=T_\nu*\mathcal{N}_{M^{\nu}} *\mathcal{N}_{M^{\mu+\mu_0}}*\w\phi^-=T_\nu *\mathcal{N}_{M^{\mu+\mu_0}}*\w\phi^-,\quad \nu=\mu, \mu+1.
$$
It follows that
$$
q_{\mu+1}(k)-q_{\mu}(k)=(T_{\mu+1}-T_{\mu})*(\mathcal{N}_{M^{\mu+\mu_0}}*\w\vp^-)(k).
$$
Then, using~\eqref{DefS}, we obtain
\begin{equation}\label{121++}
  \begin{split}
      \|\{q_{\mu+1}(k)&-q_\mu(k)\}_{k}\|_{\ell_\infty}=\|(T_{\mu+1}-T_\mu)*(\mathcal{N}_{M^{\mu+\mu_0}}*\w\vp^-)\|_{}\\
      &\le \|\mathcal{N}_{M^{\mu+\mu_0}}*\w\vp^-\|_1\|T_{\mu+1}-T_\mu\|_{}\le \a(M^{\mu+\mu_0})\|T_{\mu+1}-T_\mu\|_{}\\
      &\le c_1\big(\alpha(M^{\mu}) E_{\delta M^\mu}(f)_{} +  \alpha (M^{\mu+1})E_{\delta M^{\mu+1}}(f)_{}\big),
   \end{split}
\end{equation}
which after the corresponding summation implies~\eqref{elem1}.

Next, it is clear that there exists $\nu(\d)\in \N$ such that $E_{\delta M^\mu}(f)_{}\le E_{M^{\mu-\nu(\delta)}}(f)_{}$ and $\alpha(M^{\mu})\le C(\delta)\alpha(M^{\mu-\nu(\delta)})$
 for all big enough $\mu$.
Thus, if  $f\in \mathbb{B}_{M}^{\a(\cdot)}$, then it follows from~\eqref{elem1} that
$\{\{q_\mu(k)\}_k\}_{\mu=1}^\infty$ is a Cauchy sequence in
$\ell_\infty$.  Fortiori, for every $k\in\zd$, the sequence $\{q_\mu(k)\}_{\mu=1}^\infty$  has a limit.

Let now $\alpha=\const$ and $f\in C$. For every $\mu', \mu''\in\n$ there exists
$\nu\in\n$ such that both the functions $\h{T_{\mu'}}$ and $\h{T_{\mu''}}$ are supported
in $M^{*\nu}\T^d$, and similarly to~\eqref{121++}, we have
$$
 \|\{q_{\mu'}(k)-q_{\mu''}(k)\}_{k}\|_{\ell_\infty}\le\|(T_{\mu'}-T_{\mu''})*(\mathcal{N}_{M^\nu}*\w\vp^-)\|_{}\le c_2\big(E_{\delta M^{\mu'}}(f) +
 E_{\delta M^{\mu''}}(f)\big).
$$
Thus, again $\{\{q_\mu(k)\}_k\}_{\mu=1}^\infty$ is a Cauchy sequence in
$\ell_\infty$, and  every  sequence $\{q_\mu(k)\}_{\mu=1}^\infty$  has a limit.

Let us  check that the limit of $\{\{q_\mu(k)\}_k\}_{\mu=1}^\infty$ in $\ell_\infty$ does not depend on the choice of $T_\mu$ and $\d$.  Let $\d'\in (0,1]$ and
$T'_\mu\in \mathcal{B}_{\d' M^\mu}$  be such that
$ \|f-T'_\mu\|\le c'(d){E}_{\d' M^\mu}(f)$ and $q_{\mu}'(k)=T'_\mu*(\mathcal{N}_{M^\mu}*\w\vp^-)(k)$.
Since both the functions
$T_\mu$ and $T'_{\mu}$ belong to $\mathcal{B}_{M^{\mu}}$, repeating the
arguments of the proof of inequality~\eqref{121++} with $T'_\mu$ instead of $T_{\mu+1}$ and
$0$ instead of $\mu_0$, we obtain
\begin{equation*}
  \begin{split}     	
	\|\{q_{\mu}'(k)&-q_\mu(k)\}_{k} \|_{\ell_\infty}\le
c_3\alpha(M^\mu) \|T'_{\mu}-T_\mu\|_{}
   \le  c_4 \alpha (M^{\mu}) (E_{\delta M^\mu}(f)_{}+E_{\delta' M^\mu}(f)_{}).
   \end{split}
\end{equation*}
It follows that
$\|\{q_{\mu}'(k)-q_\mu(k)\}_{k} \|_{\ell_\infty}\to 0$
as $\mu\to\infty$, which yields the independence  on the choice of $T_\mu $ and $\delta$.~$\Diamond$

\medskip

The proof of the next lemma is obvious.

\begin{lem}
\label{lemKK1}
Let $\phi\in {\cal L}_\infty$, and  $\{a_k\}_{k\in\zd}\in\ell_\infty$.  Then
\begin{equation*}
  \bigg\Vert \sum_{k\in \Z^d} a_k \phi_{0k}\bigg\Vert_{}
\le  \|\phi\|_{\mathcal{L}_\infty} \left\Vert \{a_k\}_{k}\right\Vert_{\ell_\infty}.
\end{equation*}
\end{lem}

\begin{lem}
\label{lem99}
Let $f\in {C}$ and $\w\vp\in \mathcal{S}_{\const;M}'$ for some $M\in \mathfrak{M}$. Then
$$
\|\{\langle f, \w\phi_{0k}\rangle\}_k\|_{\ell_\infty}\le c\,\|f\|_{},
$$
where $c$ does not depend on $f$.
\end{lem}
{\bf Proof.} By Lemma~\ref{lem1}, for any $\varepsilon>0$, there exists a function   $T_\mu\in \mathcal{B}_{M^\mu}$
such that   $\|f-T_\mu\|_{}\le c_1{E}_{M^\mu}(f)_{}$ and
\begin{equation}\label{yu1}
  \|\{\langle f, \w\phi_{0k}\rangle-T_\mu*(\mathcal{N}_{M^\mu}*\w\vp^-)(k)\}_k\|_{\ell_\infty}<\varepsilon.
\end{equation}
Moreover, due to~\eqref{DefS}, we have
$$
|T_\mu*(\mathcal{N}_{M^\mu}*\w\vp^-)(k)|
\le\Vert T_\mu*(\mathcal{N}_{M^\mu}*\w\vp^-) \Vert_{} \le c_2\Vert T_{\mu}\Vert_{}
\le c_3\|f\|_{}.
$$
Combining this with~\eqref{yu1}, we prove the lemma.~$\Diamond$

\medskip

In the next two lemmas, we recall some basic inequalities for the best approximation and moduli of smoothness.

\begin{lem}\label{lemJ} {\sc (See~\cite[5.2.1 (7)]{Nik} or~\cite[5.3.3]{Timan})}
Let $f\in C$ and $s\in \N$. Then
$$
E_{I}(f)_{}\le c\,\omega_{s}(f, 1)_{},
$$
where $c$ is a constant independent of $f$.
\end{lem}

\begin{lem}\label{lemNS} {\sc (See~\cite{Wil} or~\cite{KT20})}
  Let $s\in \N$ and $T\in \mathcal{B}_{I}$. Then
\begin{equation*}
   \sum_{[\beta]=s} \Vert D^\beta T\Vert_{}\le c\,\omega_s(T,1)_{},
\end{equation*}
where $c$ does not depend on $T$.
\end{lem}

\section{Main results}
Let $M\in \mathfrak{M}$, $\a\in \A_M$,  $\w\vp$ belong to $\mathcal{S}_{\a;M}'$ and $\phi\in\mathcal{L}C$. In what follows, we  understand
$\langle f,\w\vp_{jk}\rangle$ in the sense of Definition~\ref{def0}. 	
Thus, the quasi-projection  operators
$$
Q_j(f,\phi,\w\phi)=\sum_{k\in\zd}\langle f,\w\vp_{jk}\rangle\phi_{jk}
$$
are  defined for all $f\in \mathbb{B}_{M}^{\a(\cdot)}$.
By Lemmas~\ref{lem1} and~\ref{lem99}, we have that $\{\langle f,\w\vp_{jk}\rangle\}_k\in\ell_\infty$.
This, together with Lemma~\ref{lemKK1}, implies that
the operator $Q_j(f,\phi,\w\phi)$ is well defined.

\subsection{The case of weak compatibility of  $\phi$ and $\w\phi$}

In this subsection,
we give  error estimates for the quasi-projection operators associated with weakly compatible
 $\phi$ and  $\w\phi$.

\begin{theo}
\label{corMOD1'+}
Let $M\in \mathfrak{M}$, $\a\in\A_M$,  $s\in \N$, and $\delta\in (0,1]$.
Suppose
\begin{itemize}
  \item[$1)$]  $\w\phi\in \S'_{\a;M}$ and $\phi\in\mathcal{L}C$;
  \item[$2)$] the Strang-Fix condition of order $s$  holds for $\phi$;
  \item[$3)$] $\phi$ and ${\w\phi}$ are weakly  compatible of order~$s$;
  \item[$4)$] $\eta_\delta D^{\beta}\overline{\h\phi}{\h{\w\phi}} \in W_0$ and $  \eta_\delta D^{\beta}\h{\phi}(\cdot+l)\in W_0$ for all $\beta \in \Z_+^d$, $[\beta]=s$, and  $l\in\zd\setminus \{\nul\}$;
 \item[$5)$]  $ \sum_{l\ne\nul} \|\eta_\delta D^{\beta}\h{\phi}(\cdot+l)\|_{ W_0}<\infty$ for all $\beta \in \Z_+^d$, $[\beta]=s$.
\end{itemize}
Then, for any $f\in \mathbb{B}_{M}^{\a(\cdot)}$, we have
\begin{equation}\label{110}
\begin{split}
\Vert f - Q_j(f,\phi,\w\phi) \Vert_{}\le c\(\Omega_s(f, M^{-j})_{}+\sum_{\nu=j}^\infty \alpha(M^{\nu-j}) E_{M^\nu}(f)_{}\).
\end{split}
\end{equation}
Moreover, if   $\w\vp\in \mathcal{S}_{\const;M}'$ and $f\in {C}$, then
\begin{equation}\label{110K}
  \Vert f - Q_j(f,\phi,\w\phi) \Vert_{}\le c\,\Omega_s(f,M^{-j})_{}.
\end{equation}
In the above inequalities, the constant $c$ does not depend on $f$ and $j$.
\end{theo}
{\bf Proof.} First we will prove the inequality
\be
\label{102}
\bigg\|T-\sum_{k\in\zd}\langle  T, \w\phi_{0k}\rangle\phi_{0k}\bigg\|_{} \le c_1 \sum_{[\beta]=s}\|D^\beta T\|_{},\quad T\in\mathcal{B}_{\d I}.
\ee
Set $\w\Phi(x)=T*(\mathcal{N}_\delta*\w\phi^-)(-x)$. Since
$\w\Phi\in L_\infty$ and $\phi\in \mathcal{L}C\subset L_1$, the function
$$
g_x(y)=\sum_{\nu\in \Z^d}\w\Phi(-y+\nu)\vp(x-y+\nu), \quad x\in\rd,
$$
is continuous and summable on $\td$.
Let us check that the Fourier series of $g_x$ is absolutely convergent, i.e.,
\be
\label{801}
\sum_{k\in \Z^d}|\h{g_x}(k)|=\sum_{k\in \Z^d}\Big|\int_{\R^d} \w\Phi(-y)\vp(x-y) e^{-2\pi i(y,k)}dy\Big|<\infty.
\ee
Let $k\in\zd$, $k\ne\nul$, be fixed. Denoting $\Phi_k(y)=\vp(y) e^{-2\pi i(y,k)}$, $e_k(x)=e^{2\pi i(k,x)}$, we get
\begin{equation}\label{803}
\begin{split}
    \int_{\R^d} &\w\Phi(-y)\vp(x-y) e^{-2\pi i(y,k)}dy=\int_{\R^d}T*(\mathcal{N}_\d*\w\vp^-)(y)\vp(x-y) e^{-2\pi i(y,k)}dy\\
&=e_k(x)((T*\mathcal{N}_\d)*(\mathcal{N}_\d*\w\vp^-)*\Phi_k(x)
=e_k(x)(T*(\mathcal{N}_\d*\mathcal{N}_\d*\w\vp^-)*\Phi_k)(x).
\end{split}
\end{equation}
By condition 2) and Taylor's formula, we have
\begin{equation*}
  \h{\phi} (\xi+l) = \sum_{[\beta]=s} \frac{s}{\beta!}  \xi^{\beta} \int_0^1 (1-t)^{s-1} D^{\beta}\h{\phi}( t\xi+l) d t.
\end{equation*}
Using this equality  and setting
 $K_{t,\delta,k}(x):=\mathcal{F}^{-1}\(\eta_{\d}(t\xi) D^\b {\h{\vp}}(t\xi+k)\)(x)$, we obtain
\begin{equation}
\label{i11}
  \begin{split}
     (\mathcal{N}_\d*\mathcal{N}_\d*\w\vp^-)&*\Phi_k(x)=\int_{\R^d} \eta_\d^2(\xi) \h\vp(\xi+k)\overline{\h{\w\vp}}(\xi)e^{2\pi i(\xi,x)}d\xi\\
    &=\sum_{[\b]=s}\frac{s}{\b!}\int_0^1(1-t)^{s-1}\int_{\R^d} \xi^\b\eta_\d^2(\xi)
 \overline{\h{\w\vp}}(\xi)\eta_{2\d}(t\xi) D^\b \h\vp(t\xi+k)  e^{2\pi i(\xi,x)}d\xi dt\\
     &=\sum_{[\b]=s}\frac{s}{\b!(2\pi i)^\b}\int_0^1(1-t)^{s-1}\int_{\R^d} \h{D^\b \mathcal{N}_{\d}} (\xi)\eta_\d(\xi)\overline{\h{\w\vp}}(\xi) \h{K_{t,2\delta,k}}(\xi)e^{2\pi i(\xi,x)}d\xi dt\\
     &=\sum_{[\b]=s}\frac{s}{\b!(2\pi i)^\b}\int_0^1(1-t)^{s-1}(D^\b \mathcal{N}_{\d}*(N_\d*\w\vp^-)*K_{t,2\delta,k})(x)dt.\\
  \end{split}
\end{equation}
Since, the function $\mathcal{N}_\d*\w\vp^-$ is summable, it follows that
\begin{equation}\label{i12}
  \begin{split}
\Vert T*(\mathcal{N}_\d*\mathcal{N}_\d*\w\vp^-)*\Phi_k\Vert_{}&
\le c_2 \sum_{[\b]=s}\int_0^1(1-t)^{s-1}\Vert T*D^\b \mathcal{N}_{\d}*(N_\d*\w\vp^-)*K_{t,2\d,k}\Vert_{} dt
\\
&\le c_3
\sum_{[\b]=s} \sup_{t\in (0,1)}\Vert T*D^\b \mathcal{N}_{\d}*K_{t,2\d,k}\Vert_{}.
  \end{split}
\end{equation}
Due to condition 4), we derive
\begin{equation}\label{i12+}
  \begin{split}
\Vert T*D^\b \mathcal{N}_{\d}*K_{t,2\d,k}\Vert_{}
&=\Vert D^\b T* \mathcal{N}_{\d}*K_{t,2\d,k}\Vert_{}\le c_4\Vert D^\b T*K_{t,2\d,k}\Vert_{}
\\
&= c_4\Vert D^\b T*\mathcal{N}_{t^{-1}\d}*K_{t,\d,k}\Vert_{}\le c_4\Vert D^\b T*K_{t,\d,k}\Vert_{}
\\
&\le c_4\Vert K_{t,\d,k}\Vert_1\Vert D^\b T\Vert_{}\le c_5 \|\eta_{\d} D^{\beta}\h{\phi}(\cdot+k)\|_{ W_0}
 \Vert D^\b T\Vert_{}.
   \end{split}
\end{equation}
Thus, inequalities \eqref{i12}, \eqref{i12+}, and  condition 5) yield
\be
\label{802}
\Big\|\sum_{k\neq \nul}e_k(T*(\mathcal{N}_\d*\mathcal{N}_\d*\w\vp^-)*\Phi_k)\Big\|_{}\le c_6 \sum_{[\b]=s} \Vert D^\b T\Vert_{},
\ee
which completes the proof of~\eqref{801}. Therefore, by the Poisson summation formula  and~\eqref{803}, we have
\begin{equation}\label{804}
  \begin{split}
    \sum_{k\in\Z^d} &(T*(\mathcal{N}_\d*\w\vp^-))(k)\phi_{0k}(x)=\sum_{k\in \Z^d}\w\Phi(k)\vp(x+k)\\
&=\sum_{k\in\zd}\int_{\R^d} \w\Phi(-y)\vp(x-y) e^{-2\pi i(y,k)}dy=\sum_{k\in\zd}e_k(T*(\mathcal{N}_\d*\mathcal{N}_\d*\w\vp^-)*\Phi_k),
  \end{split}
\end{equation}
and hence
\begin{equation}\label{i8}
  \begin{split}
T-\sum_{k\in\Z^d} (T*(\mathcal{N}_\d*\w\vp^-))(k)\phi_{0k}=T*&(\mathcal{N}_\d-(\mathcal{N}_\d*\w\vp^-)*\vp)\\
&+\sum_{k\neq \nul}e_k(T*(\mathcal{N}_\d*\mathcal{N}_\d*\w\vp^-)*\Phi_k)=:I_1+I_2.
  \end{split}
\end{equation}
The term $I_2$ is already estimated in~\eqref{802}.
Consider the term $I_1$. By condition 4) and  Taylor's formula, there holds
\begin{equation*}
  \h{\phi} (\xi)\overline{\h{\w\phi}(\xi)} = 1 + \sum_{[\beta]=s}
	\frac{s}{\beta!}  \xi^{\beta}
	\int_0^1 (1-t)^{s-1} D^{\beta}\h{\phi}\overline{\h{\w\phi}}( t \xi) d t.
\end{equation*}
Using this equality, we obtain
\begin{equation}\label{i9}
  \begin{split}
     \mathcal{N}_\d(x)&-((\mathcal{N}_\d*\w\vp^-)*\vp)(x)=\int_{\R^d} \eta_\d(\xi) (1-\h\vp(\xi)\overline{\h{\w\vp}}(\xi))e^{2\pi i(\xi,x)}d\xi\\
     &=\sum_{[\b]=s}\frac{s}{\b!}\int_0^1(1-t)^{s-1}\int_{\R^d} \xi^\b\eta_{\d}(\xi) \eta_{2\d}(t\xi)D^\b \h\vp\overline{\h{\w\vp}}(t\xi)e^{2\pi i(\xi,x)}d\xi dt\\
     &=\sum_{[\b]=s}\frac{s}{\b!(2\pi i)^\b}\int_0^1(1-t)^{s-1}\int_{\R^d} \h{D^\b \mathcal{N}_{\d}} (\xi) \h{K_{t,2\d}}(\xi)e^{2\pi i(\xi,x)}d\xi dt\\
     &=\sum_{[\b]=s}\frac{s}{\b!(2\pi i)^\b}\int_0^1(1-t)^{s-1}(D^\b \mathcal{N}_{\d}*K_{t, 2\d})(x)dt,\\
  \end{split}
\end{equation}
where $K_{t,\d}(x)=\mathcal{F}^{-1}\(\eta_{\d}(t\xi) D^\b \h\vp\overline{\h{\w\vp}}(t\xi)\)(x)$.
Next, using  condition 4) and the same arguments as in~\eqref{i12+}, we derive
\begin{equation}\label{i10}
  \begin{split}
\Vert I_1 \Vert_{}&\le c_{7} \sum_{[\b]=s}\int_0^1(1-t)^{s-1}\Vert T*D^{\b} \mathcal{N}_\d*K_{t,2\d}\Vert_{} dt
\\
&\le c_{8} \sum_{[\b]=s} \sup_{t\in (0,1)}\Vert D^\b T* K_{t,\d} \Vert_{}\le c_{9}\sum_{[\b]=s} \Vert D^\b T\Vert_{}.
  \end{split}
\end{equation}
Combining this with~\eqref{i8} and~\eqref{802}, we obtain~\eqref{102}.

Now let the functions $T_\mu$, $\mu \in \Z_+$, be as in Definition~\ref{def0} and $T=T_0$.
It follows from \eqref{102} and
Lemma~\ref{lemNS}  that
\be
\label{102+}
\bigg\|T-\sum_{k\in\zd}\langle  T, \w\phi_{0k}\rangle\phi_{0k}\bigg\|_{} \le c_{10} \omega_s(T,1)\le c_{11}\(\omega_s(f,1)+E_{\d I}(f)\).
\ee
If $f\in \mathbb{B}_{M}^{\a(\cdot)}$, then, by  Definition~\ref{def0}, using
Lemmas~\ref{lem1} and~\ref{lemKK1},  we derive
\begin{equation}\label{103}
  \begin{split}
     \|Q_0(f-T,\vp,\w\vp)\|_{}&\le c_{12}\|\{\langle f-T, \w\phi_{0k}\rangle\}_k\|_{\ell_\infty}\le c_{13}\sum_{\mu=0}^\infty\|\{\langle T_{\mu+1}-T_\mu, \w\phi_{0k}\rangle\}_k\|_{\ell_\infty}\\
&\le c_{14}\sum_{\mu=\nu}^\infty \alpha (M^\mu) E_{\d M^\mu}(f)_{}.
  \end{split}
\end{equation}
Combining~\eqref{102+}  and~\eqref{103} with the inequality
\begin{equation}\label{_1}
  \begin{split}
     \|f-Q_0(f,\vp,\w\vp)\|\le \|T-Q_0(T,\vp,\w\vp)\|+\|Q_0(f-T,\vp,\w\vp)\|+\|f-T\|
  \end{split}
\end{equation}
and taking into account that $\|f-T\|\le c(d) E_{\delta I}(f)$, we obtain
\be
\label{104}
\Vert f - Q_0(f,\phi,\w\phi) \Vert_{}\le c_{15} \(\omega_s(f,1)_{}+\sum_{\nu=0}^\infty \alpha(M^{\nu}) E_{\d M^\nu}(f)_{}\).
\ee
Since there exists $\nu_0=\nu(\d)\in \N$ such that $E_{\d M^\nu}(f)_{}\le E_{M^{\nu-\nu_0}}(f)$
and $\alpha(M^{\nu})\le c(\d)\alpha(M^{\nu-\nu_0})$
for all  $\nu>\nu_0$, applying Lemma~\ref{lemJ} and the inequality $\omega_s(f,\l)_{}\le (1+\l)^s \omega_s(f,1)_{}$ (see, e.g.,~\cite{KT20}) to the first~$\nu_0$ terms of the sum, we get~\eqref{110} for $j=0$.

If $\w\vp\in \mathcal{S}_{\const;M}'$ and $f\in {C}$, then, by Definition~\ref{def0} and Lemma~\ref{lem99},  we have
\be
\label{105}
 \|Q_0(f-T,\vp,\w\vp)\|\le c_{16}\|\{\langle f-T, \w\phi_{0k}\rangle\}_k\|_{\ell_\infty}\le
c_{17}E_{\delta I}(f).
\ee
Combining this with~\eqref{102+} and~\eqref{_1},  taking into account that $\|f-T\|\le c(d) E_{\delta I}(f)$, using Lemma~\ref{lemJ} and the properties of moduli of smoothness, we get~\eqref{110K} for $j=0$.

Thus, our theorem is proved for the case $j=0$. To prove~\eqref{110} and~\eqref{110K}
for arbitrary $j$, it remains to note that
$$
\bigg\|f-\sum\limits_{k\in\zd}\langle f,\widetilde\phi_{jk}\rangle \phi_{jk}\bigg\|_{}=
\bigg\|f(M^{-j}\cdot)-\sum\limits_{k\in\zd}\langle f(M^{-j}\cdot),\widetilde\phi_{0k}\rangle \phi_{0k}\bigg\|,
$$
$$
E_{ M^\nu}(f(M^{-j}\cdot))_{}=E_{ M^{\nu+j}}(f),
$$
$$
\omega_s(f(M^{-j}\cdot), 1)_{}=\Omega_s(f, M^{-j}),
$$
and that $f(M^{-j}\cdot)\in \mathbb{B}_{M}^{\a(\cdot)}$
whenever $f\in \mathbb{B}_{M}^{\a(\cdot)}$.~~$\Diamond$

\begin{coro}
\label{coro1NN}
Let $s\in \N$, $\delta\in (0,1]$, $M\in \mathfrak{M}$, $\a\in \A_M$, $\w\phi\in \S'_{\a;M}$ and
$\phi\in {\cal L}C$.
Suppose that conditions 2) and 3) of Theorem~\ref{corMOD1'+} are satisfied and, additionally,
for some $k\in \N$, $k>d/2$,
$$
\overline{\h{\phi}}{\h{\w\vp}}\in W_2^{s+k}(2\delta \T^d),
\quad \h\phi(\cdot+l) \in W_2^{s+k}(2\delta \T^d)\quad\text{for all}\quad l\in\Z^d\setminus \{{\bf 0}\},
$$
and
$$
\sum_{l\neq \bf{0}}\Vert D^\beta \h \vp(\cdot+l)\Vert_{L_2(2\d \T^d)}^{1-\frac{d}{2k}}<\infty\quad \text{for all}\quad \beta \in \Z_+^d,\quad [\beta]=s.
$$
Then inequalities~\eqref{110} and~\eqref{110K} hold true.
\end{coro}

The proof of Corollary~\ref{coro1NN} easily follows from Theorem~\ref{corMOD1'+} and the Beurling-type sufficient condition for belonging to Wiener's algebra given in~\cite[Theorem~6.1]{LST}. Note that more general efficient sufficient conditions on smoothness of~$\h\vp$ and $\h{\w\vp}$ can be obtained by exploiting the results of the papers~\cite{Kol}, \cite{KL13}, \cite{KL17}.

\bigskip

\noindent\textbf{Example 1.} Let $d=2$,  $s=2$, $f\in {C}$, and
$ Q_j(f,\phi,\w\phi)$   be a  mixed sampling-Kantorovich  quasi-projection operator associated with
$\phi(x)=4^{-d}\prod\limits_{l=1}^{d}\sinc^3(x_l/4)$ and
$$
\w\phi(x)=\prod\limits_{l=1}^{d'}\chi_{[-1/2, 1/2]}(x_l)\prod\limits_{l=d'+1}^{d}\delta(x_l).
$$
Thus,
$$
Q_j(f,\vp,\w\vp)(x)=\sum_{k\in\zd}\int\limits_{k_1-1/2}^{k_1+1/2}dt_1\dots \int\limits_{k_{d'}-1/2}^{k_{d'}+1/2}dt_{d'}
f(M^{-j}(t+k))\Big|_{t_{d'+1}=\dots=t_d=0}\, \phi (Mx-k).
$$
It is easy to see that all assumptions of Theorem~\ref{corMOD1'+} for the case
$\w\vp\in \S_{{\rm const};M}'$ are satisfied, which implies
$$
\|f-Q_j(f,\vp,\w\vp)\|_{}\le c\,\Omega_2(f,M^{-j})_{}.
$$

\medskip

In Theorem~\ref{corMOD1'+} above, the uniform error estimates for the operators $Q_j(f,\vp,\w\vp)$ are given in terms classical  moduli of smoothness. Similar estimates in $L_p$-norm, $p<\infty$, were obtained in our recent paper~\cite{KS20+}, where we also derived the corresponding lower estimates of approximation in terms of the same moduli of smoothness. It turns out that in the uniform metric the classical moduli of smoothness are not suitable to obtain sharp estimates of approximation (or two-sided inequalities) in the general case (see, e.g,~\cite[Ch.~9]{TB}). However, this situation can be improved by using some special measures of smoothness (e.g. $K$-functionals and their realizations as well as some special moduli of smoothness) instead of the classical moduli of smoothness (see, e.g.,~\cite[Ch.~9]{TB}, \cite{ARS20}, \cite{KTmemo}). For our purposes, it is convenient to use realizations of $K$-functionals.

To give the definition of the realization of $K$-functionals,  we need some additional notation.
We say that a function $\rho$ belongs to the class $\mathcal{H}_s$, $s>0$,
if  $\rho\in C^\infty(\R^d\setminus \{\bf{0}\})$ is a homogeneous function of degree~$s$, i.e.,
$$
\rho(\tau \xi)=\tau^s \rho(\xi),\quad \xi\in\R^d.
$$

Any function $\rho\in \mathcal{H}_s$ generates the Weyl-type differentiation operator as follows:
$$
\mathcal{D}(\rho)g:=\F^{-1}\(\rho \h g\),\quad g\in \S.
$$
For functions $T\in \B_{M^j}$, we define $\mathcal{D}(\rho)T$ by
$$
\mathcal{D}(\rho)T:=(\mathcal{D}(\rho)\mathcal{N}_{M^j})*T.
$$
Note that the operator is well defined because the function $\F^{-1}(\mathcal{D}(\rho)\mathcal{N}_{M^j})$ is summable on $\R^d$,  see, e.g.,~\cite{RS}.

Let us give several important  examples of the Weyl-type  operators:

\begin{enumerate}
  \item[{(1)}]\; the linear differential operator
$$
P_m({D})f=
\sum_{[k]=m}a_k \frac{\partial^{k_1+\cdots+k_d}}{\partial x_1^{k_1}
\cdots \partial x_d^{k_d}}f,
$$
which corresponds to $\rho(\xi)=\sum_{[k]=m}a_k (i\xi_1)^{k_1}\dots(i\xi_d)^{k_d}$;
\item [{(2)}]\; the fractional Laplacian  $(-\Delta)^{s/2}f$ (here $\rho(\xi)=|\xi|^s$, $\xi\in\R^d$);
\item [{(3)}] \;
  the classical
Weyl derivative $f^{(s)}$ (here $\rho(\xi)=(i \xi)^s$, $\xi\in\R$).
\end{enumerate}

The realization of the $K$-functional generated by the function $\rho\in \mathcal{H}_s$ is defined by
$$
\mathcal{R}_\rho(f,M^{-j})_{}=\inf_{T\in \B_{M^j}}\{\|f-T\|_{}+\|\mathcal{D}(\rho(M^{*-j}\cdot))T\|_{}\}.
$$
In many cases, the realization $\mathcal{R}_\rho(f,M^{-j})_{}$ is equivalent to the corresponding modulus of smoothness or $K$-functional (see, e.g.,~\cite{ARS20}, \cite{KT20}, \cite{KTmemo}).
Let us mention some results in this direction. For example, if $d=1$ and $\rho(\xi)=(i\xi)^s$, then
$$
\mathcal{R}_\rho(f,M^{-j})_{}\asymp \omega_s(f,M^{-j})_{},
$$
where $\omega_s(f,M^{-j})_{}$ is the standard fractional modulus of smoothness defined in~\eqref{fm}.
If $d=1$ and $\rho(\xi)=|\xi|^s$, then
$$
\mathcal{R}_\rho(f,M^{-j})_{}\asymp \omega_s^\Delta(f,M^{-j})_{},
$$
where
$$
\omega_s^\Delta(f,M^{-j})_{}=\sup_{|M^j h|\le 1} \bigg\|\sum_{\nu\neq 0}\frac{f(\cdot+\nu h)-f(\cdot)}{|\nu|^{s+1}}\bigg\|_{}.
$$
In the case $d\ge 1$ and $\rho(\xi)=|\xi|^2$, we have
$$
\mathcal{R}_\rho(f,M^{-j})_{}\asymp \w\Omega_2(f,M^{-j})_{},
$$
where
$$
\w\Omega_2(f,M^{-j})_{}=\sup_{|M^j h|\le 1} \bigg\|\sum_{j=1}^d\(f(\cdot+h{\rm e}_j)-2f(\cdot)+f(\cdot-h{\rm e}_j)\)\bigg\|_{}
$$
and $\{{\rm e}_j\}_{j=1}^d$ is the standard basis in $\R^d$.

\begin{theo}
\label{corMOD1+f++}
Let $s>0$, $\rho\in \mathcal{H}_s$, $\delta\in (0,1/2)$, $M\in \mathfrak{M}$, and  $\a\in \A_M$.
Suppose that $\w\phi\in \S'_{\a;M}$ and  $\phi\in\mathcal{L}C$ satisfy the following conditions:
\begin{itemize}
  \item[$1)$] $\supp\h\vp \subset\T^d$;
  \item[$2)$] $\eta_\delta\frac{1-\overline{\h\phi}{\h{\w\phi}}}{\rho} \in W_0$.
\end{itemize}
Then, for any $f\in \mathbb{B}_{M}^{\a(\cdot)}$, we have
\begin{equation*}
\begin{split}
\Vert f - Q_j(f,\phi,\w\phi) \Vert_{}\le c \(\mathcal{R}_\rho(f, M^{-j})_{}+\sum_{\nu=j}^\infty \a(M^{\nu-j}) E_{M^\nu}(f)_{}\).
\end{split}
\end{equation*}
Moreover, if   $\w\vp\in \mathcal{S}_{\const;M}'$ and $f\in {C}$, then 
\begin{equation*}
  \Vert f - Q_j(f,\phi,\w\phi) \Vert_{}\le c\, \mathcal{R}_\rho(f,M^{-j})_{}.
\end{equation*}
In the above inequalities, the constant $c$ does not depend on $f$ and $j$.
\end{theo}
\textbf{Proof.}
Analyzing the arguments of the proof of Theorem~\ref{corMOD1'+}, one can see that it suffices to verify that
\begin{equation}\label{eqNN1}
  \|T*\(\mathcal{N}_\d*\mathcal{N}_\d-(\mathcal{N}_\d*\mathcal{N}_\d*\w\vp)*\vp\)\|_{}\le c_1\|\mathcal{D}(\rho)T\|_{}
\end{equation}
for any $T \in \mathcal{B}_{\d I}$ such that $\Vert f-T\Vert_{}\le c(d)E_{\d I}(f)_{}$.

We have
\begin{equation*}
  \begin{split}
      \mathcal{N}_\d*\mathcal{N}_\d(x)-(\mathcal{N}_\d*\mathcal{N}_\d*\w\vp)*\vp(x)&=\int_{\R^d}\eta_\d^2(\xi)\(1-\h\vp(\xi)\overline{\h{\w\vp}}(\xi)\)e^{2\pi i(\xi,x)}d\xi\\
      &=\int_{\R^d}\eta_\d(\xi)\frac{1-\h\vp(\xi)\overline{\h{\w\vp}}(\xi)}{\rho(\xi)}\rho(\xi)\eta_\d(\xi)e^{2\pi i(\xi,x)}d\xi\\
      &=(\mathcal{D}(\rho)\mathcal{N}_\d)*K_\d(x),
   \end{split}
\end{equation*}
where $K_\d(x)=\F^{-1}\(\eta_\d\frac{1-\h\vp(\xi)\overline{\h{\w\vp}}(\xi)}{\rho(\xi)}\)(x)$. Thus, using condition 2), we derive
\begin{equation*}
  \begin{split}
     \|T*\(\mathcal{N}_\d*\mathcal{N}_\d-(\mathcal{N}_\d*\mathcal{N}_\d*\w\vp)*\vp\)\|_{}=\|T*(\mathcal{D}(\rho)\mathcal{N}_\d)*K_\d\|_{}\le c_2\|T*(\mathcal{D}(\rho)\mathcal{N}_\d)\|_{}=c_2\|\mathcal{D}(\rho)T\|_{},
   \end{split}
\end{equation*}
which proves the theorem.~~$\Diamond$

\bigskip

In the next theorem, we obtain lower estimates of the approximation error by the quasi-projection operators $Q_j(f,\phi,\w\phi)$. Note that such type of estimates are also called strong converse inequalities, see, e.g.,~\cite{DI}.

\begin{theo}\label{corMOD2conv}
Let $s>0$, $\rho\in \mathcal{H}_s$, $M\in \mathfrak{M}$,  and $\a\in \A_M$.
Suppose that $\w\phi\in \S'_{\a;M}$ and  $\phi\in\mathcal{L}C$ satisfy the following conditions:
\begin{itemize}
  \item[$1)$] $\supp \h\vp\subset \T^d$;
  \item[$2)$]  $\eta\frac{\rho}{1-\overline{\h\phi}{\h{\w\phi}}} \in W_0$.
\end{itemize}
Then, for any $f\in \mathbb{B}_{M}^{\a(\cdot)}$, we have
\begin{equation}\label{110inv}
\begin{split}
\mathcal{R}_\rho(f, M^{-j})_{}\le c\,\Vert f - Q_j(f,\phi,\w\phi) \Vert_{}+c\sum_{\nu=j}^\infty \alpha(M^{\nu-j}) E_{M^\nu}(f)_{}.
\end{split}
\end{equation}
Moreover, if $\w\vp\in \mathcal{S}_{\const;M}'$ and $f\in {C}$, then
\begin{equation}\label{110Kinv}
  \mathcal{R}_\rho(f,M^{-j})_{}\le c\,\Vert f - Q_j(f,\phi,\w\phi) \Vert_{}.
\end{equation}
In the above inequalities, the constant $c$ does not depend on $f$ and $j$.
\end{theo}
\textbf{Proof.} As in the proof of the previous theorems, it suffices to consider only the case $j=0$.
Let $T \in \mathcal{B}_{I}$ be such that $\Vert f-T\Vert_{}\le c(d)E_{I}(f)_{}$.
Due to the same arguments as in the proof of Theorem~\ref{corMOD1'+}, we have~\eqref{i8}, which takes now the following form
$$
T-\sum\limits_{k\in \Z^d}
\langle T, \w\phi_{0k} \rangle \phi_{0k}=T-\sum_{k\in\Z^d} (T*(\mathcal{N}_1*\w\vp^-))(k)\phi_{0k}=T*(\mathcal{N}_1-(\mathcal{N}_1*\w\vp^-)*\vp).
$$
Thus, denoting $K=\F^{-1}(\eta\frac{\rho}{1-\overline{\h\phi}{\h{\w\phi}}})$ and using condition~2), we derive
\begin{equation*}
  \begin{split}
     \|\mathcal{D}(\rho)T\|_{}&=\|T*\mathcal{D}(\rho)\mathcal{N}_1\|_{}=\|T*\F^{-1}(\rho\eta)\|_{}
     =\bigg\|T*\F^{-1}\bigg(\eta\frac{\rho}{1-\overline{\h\phi}{\h{\w\phi}}}(1-\overline{\h\phi}{\h{\w\phi}})\eta\bigg)\bigg\|_{}\\
     &=\|T*(\mathcal{N}_1-(\mathcal{N}_1*\w\vp^-)*\vp)*K\|_{}\le c_1\|T*(\mathcal{N}_1-(\mathcal{N}_1*\w\vp^-)*\vp)\|_{}\\
     &=c_1\|T-Q_0(T,\vp,\w\vp)\|_{}.
   \end{split}
\end{equation*}
Now, by the definition of the realization, we get
\begin{equation*}
  \begin{split}
      \mathcal{R}_\rho(f,I)_{}
			&\le \Vert \mathcal{D}(\rho)T\Vert_p+E_I(f)_{}
\\
&\le c_1\(\Vert T-Q_0(T,\vp,\w\vp)\Vert_{}+E_{I}(f)_{}\)
\\
&\le c_1\(\Vert f-Q_0(f,\vp,\w\vp)\Vert_{} + \|T-f\|_{}+\Vert Q_0(T-f,\vp,\w\vp)\Vert_{}+E_{I}(f)_{}\)
\\
&\le c_2\(\Vert f-Q_0(f,\vp,\w\vp)\Vert_{}+E_{I}(f)_p+\Vert Q_0(f-T,\vp,\w\vp)\Vert_{}\)\\
&\le c_3\(\Vert f-Q_0(f,\vp,\w\vp)\Vert_{}+E_{I}(f)_{}\),
   \end{split}
\end{equation*}
where the last inequality follows from Lemmas~\ref{lemKK1} and~\ref{lem99}.
Thus, to prove \eqref{110Kinv}, it remains to note that in view of the inclusion  $\supp\mathcal{F}\(Q_0(f,\vp,\w\vp)\)\subset \supp\h\phi\subset \td$, we have
$
E_{I}(f)_{}\le \|f-Q_0(f,\vp,\w\vp)\|_{}.
$

Similarly, using~\eqref{103},
one can prove~\eqref{110inv}.~~$\Diamond$

\medskip

\begin{rem}
 Note that the conditions on  $\vp$ and $\w\vp$ in Theorems~\ref{corMOD1+f++} and~\ref{corMOD2conv} can be given also in terms of smoothness of $\h\vp$ and $\h{\w\vp}$, similarly to those given in Corollary~\ref{coro1NN}. For this, one can use the sufficient conditions for belonging to Wiener's algebra given in~\cite{LST} (see also~\cite{Kol}, \cite{KL17}).
\end{rem}

\medskip

\noindent\textbf{Example 2.}
Let $d=1$, $\vp(\xi)=\sinc^4(\xi)$, and $\w\vp(x)=\chi_{\T}(\xi)$ (the characteristic function of $\T$). Then all conditions of Theorems~\ref{corMOD1+f++}  and~\ref{corMOD2conv}  are satisfied. Therefore, for any  $f\in C$, we have
\begin{equation}\label{KKM}
 \bigg\Vert f- \sum_{k\in \Z} {M^j}\bigg(\int_{M^{-j}\T} f(M^{-j}k-t) dt\bigg) \sinc^4(M^j \cdot - k)\bigg\Vert\asymp\omega_2(f,M^{-j}),
\end{equation}
where $\asymp$ is a two-sided inequality with constants independent of $f$ and $j$.

\medskip

Now we consider approximation by the quasi-projection operators generated by the Bochner-Riesz kernel of fractional order.

\medskip

\noindent\textbf{Example 3.}
Let  $\vp(x)=R_{s}^\gamma(x):=\mathcal{F}^{-1}\((1-|3\xi|^s)_+^\g\)(x)$, $s>0$, $\g>\frac{d-1}{2}$, and $\rho(\xi)=|\xi|^s$.
  \begin{itemize}
    \item[1)] If $\w\vp(x)=\d(x)$, then for any $f\in C$, we have
\begin{equation}\label{ee61}
\begin{split}
\bigg\Vert f- {m^j}\sum_{k\in \Z^d} &f(M^{-j}k){R}_{s}^\gamma (M^j \cdot - k)\bigg\Vert_{}\asymp \mathcal{R}_\rho(f,M^{-j})_{}.
\end{split}
\end{equation}
    \item[2)] If $\w\vp(x)=\chi_{\T^d}(x)$, then  for any $f\in  C$ and $s\in (0,2]$, we have
\begin{equation}\label{ee62}
\begin{split}
\bigg\Vert f- {m^j}\sum_{k\in \Z^d} \bigg(\int_{M^{-j}\T^d} f(M^{-j}k-t) dt\bigg) &{R}_{s}^\gamma (M^j \cdot - k)\bigg\Vert_{}\asymp \mathcal{R}_\rho(f,M^{-j})_{}.
\end{split}
\end{equation}
\end{itemize}
In the above two relations, $\asymp$ is a two-sided inequality with positive constants independent of $f$ and~$j$.

\medskip

{\bf Proof.}
Estimate~\eqref{ee61} follows from Theorems~\ref{corMOD1+f++}  and~\ref{corMOD2conv} and the fact that
\begin{equation}\label{ee63+}
  \frac{\eta(\xi)(1-(1-|3\xi|^s)_+^\g)}{|\xi|^s}\in W_0\quad\text{and}\quad \frac{\eta(\xi)|\xi|^s}{1-(1-|3\xi|^s)_+^\g}\in W_0.
\end{equation}
The proof of relations~\eqref{ee63+} can be found, e.g., in~\cite{RS10}. The proof of~\eqref{ee62} is similar. In this case, instead of~\eqref{ee63+}, we use the following relations:
$$
\frac{\eta(\xi)(1-\sinc(\xi)(1-|3\xi|^s)_+^\g)}{|\xi|^s}\in W_0\quad\text{and}\quad \frac{\eta(\xi)|\xi|^s}{1-\sinc(x)(1-|3\xi|^s)_+^\g}\in W_0,
$$
which can be verified using the same arguments as in~\cite{RS10}.~~$\Diamond$

\subsection{The case of strict compatibility for  $\phi$ and $\w\phi$}

\begin{theo}\label{cor1}
Let $\d\in (0,1]$, $M\in \mathfrak{M}$, and $\a\in \A_M$.
Suppose that $\w\phi\in \mathcal{S}'_{\a;M}$ and  $\phi\in\mathcal{L}C$ satisfy the following conditions:
\begin{itemize}
  \item[$1)$] $\supp \h\vp\subset \T^d$;
  \item[$2)$] $\phi$ and ${\w\phi}$ are strictly  compatible with respect to $\d$.
\end{itemize}
Then, for any $f\in \mathbb{B}_{M}^{\a(\cdot)}$, we have
\begin{equation}
\begin{split}
\label{703}
     \Vert f - Q_j(f,\phi,\w\phi) \Vert_{}\le c\sum_{\nu=j}^\infty \a(M^{\nu-j})
E_{\d M^\nu}(f).
\end{split}
\end{equation}
Moreover, if  $\w\vp\in \S_{\const;M}'$ and $f\in {C}$,  then

\begin{equation}\label{KS000KLLL}
  \Vert f - Q_j(f,\phi,\w\phi) \Vert_{}\le c\, E_{\d M^j}(f)_{}.
\end{equation}
In the above inequalities, the constant $c$ does not depend on $f$ and $j$.
\end{theo}
{\bf Proof.}
As above, it suffices to consider only the case $j=0$. Repeating the arguments of the proof of Theorem~\ref{corMOD1'+}, we obtain
from~\eqref{i8} that
$$
T-\sum\limits_{k\in\zd}\langle T,\widetilde\phi_{0k}\rangle \phi_{0k}=0,
$$
Thus, using~\eqref{_1}, \eqref{103},  and~\eqref{105}, we prove both the statements.~~$\Diamond$

\begin{theo}\label{kot}
	Let $M\in \mathfrak{M}$ be such that $\td\subset M^*\td$, $\d\in(0,1]$, and let  $\phi$, $\w\phi$ be as in Theorem~\ref{cor1}.
If $f\in \mathcal{B}_{\d M^j}$, then
$$
f(x)=\sum_{k\in\zd}( f(M^{-j}\cdot)*\mathcal{N}_{1}*\w\vp^-)(k)\phi(M^jx-k)\quad\text{for all}\quad x\in\rd.
$$
\end{theo}
{\bf Proof.} Since obviously $E_{\d M^\nu}(f)=0$ for any $\nu\ge j$, and both the functions $f$ and $Q_j(f,\phi,\w\phi) $ are continuous, it follows from Theorem~\ref{cor1} that
 $f=Q_j(f,\phi,\w\phi)$ at each point. Thus, by Definition~\ref{def0},
 for every $x\in\rd$ we have
$$
f(x)=m^{-j/2}\sum_{k\in\zd}\langle f(M^{-j}\cdot), \w\phi_{0k}\rangle \phi_{jk}=
m^{-j/2}\sum_{k\in\zd} \lim_{\mu\to\infty} T_\mu*(\mathcal{N}_{M^{\mu}}*\w\vp^-)(k)\phi_{jk},
$$
where  $T_\mu\in \mathcal{B}_{\d M^\mu}$ is such that
$$
  \|f(M^{-j}\cdot)-T_\mu\|\le c(d){E}_{\d M^\mu}f(M^{-j}\cdot).
$$
It remains to note that for sufficiently large $\mu$ we have
$$
T_\mu*(\mathcal{N}_{M^{\mu}}*\w\vp^-)=
f(M^{-j}\cdot)*(\mathcal{N}_{M^{\mu}}*\w\vp^-)=
f(M^{-j}\cdot)*(\mathcal{N}_{1}*\w\vp^-).\quad \Diamond
$$

\noindent {\bf Acknowledgments} This research was supported by Volkswagen Foundation in framework
of the project "From Modeling and Analysis to Approximation".
The first author was partially supported by the DFG project KO 5804/1-1.

\end{document}